\documentclass[12pt]{amsart}
\usepackage{epsfig}
\usepackage{amsmath}
\usepackage{amsfonts}
\usepackage{amssymb}
\usepackage{epigraph, fancyhdr}

\def\I{\mathcal I}

\def\1{\mathbf 1}

\def\QQ{\mathbb Q}
\def\ZZ{\mathbb Z}
\def\CC{\mathbb C}

\def\hat{\widehat}
\def\tilde{\widetilde}
\def\p{\partial}
\def\a{\alpha}

\def\la{\lambda}
\def\gL{\Lambda}

\def\m{{\mathfrak m}}

\def\ch{\operatorname{ch}}

\renewcommand{\Delta}{\triangle}

\title[Mirrors]
      {Permutation-equivariant \\ quantum K-theory VI. \\
      Mirors}

\author[A. Givental]{Alexander GIVENTAL}
\thanks{This material is based upon work supported by the National 
Science Foundation under Grant DMS-1007164, and by the IBS Center for Geometry 
and Physics, POSTECH, Korea.} 

\date{July 24, 2015}

\begin{document}

\begin{abstract}

  We present here the K-theoretic version of mirror models of toric manifold.

  First, we recall the construction \cite{GiZ} of mirrors for toric manifolds, i.e. representations of toric hypergeometric functions by complex oscillating integrals. Then we repeat the construction in the K-theoretic situation and obtain complex oscillating integrals representing $q$-hypergeometric functions of toric manifolds, bundles, and super-bundles. Finally, we examine the Lagrangian varieties parameterized by critical points of the phase functions of the oscillating integrals. 

\end{abstract}

\maketitle

\section*{Hori-Vafa mirrors}

According to the Berry Principle\footnote{Formulated by Vladimir Arnold and named so after Michael Berry, known for ``the Arnold Principle'', according to which the Berry Principle is self-consistent.}, it is a sin to name inventions after their authors. To our defense, Hori--Vafa mirrors \cite{HV} were introduced in \cite{GiZ}, and here is how.

To a toric manifold $X=\CC^N//_{\omega}T^K$ one can associate the $H^*(X,\QQ)$-valued hypergeometric series
\[ I^H_X = \sum_{d\in \ZZ^K} Q^d \prod_{j=1}^N \frac{\prod_{r=-\infty}^0(u_j(p)-rz)} {\prod_{r=-\infty}^{D_j(d)}(u_j(p)-rz)},\]
where $u_j(p)=p_1m_{1j}+\cdots+p_Km_{Kj}$, $j=1,\dots, N$ are cohomology classes of toric divisors expressed in a basis $\{ p_1,\dots, p_K \}$ of $H^2(X,\ZZ)$,
and $Q^d=Q_1^{d_1}\cdots Q_K^{d_K}$ are Novikov's monomials representing degrees
$d$ of holomorphic spheres in $X$ in the dual basis of $H_2(X,\ZZ)$. The series
satisfies a certain system of $K$ differential equations, $i=1,\dots, K$:
\begin{align*}  \prod_{j: m_{ij}>0} & \prod_{r=0}^{m_{ij}-1}(m_{ij}(p_i-zQ_i\p_{Q_i})+rz) \ I_X^H = \\ 
  Q_i \prod_{j: m_{ij}<0} & \prod_{r=0}^{-m_{ij}-1}(m_{ij}(p_i-zQ_i\p_{Q_i})+rz)\ I^H_X.
  \end{align*}
As it is explained in \cite{GiZ, GiH} (and in Part V), both the series $I_X^H$ and the $D$-module it generates are by-products of a heuristic approach to $S^1$-equivariant Floer theory. For toric manifolds $X$ obtained by factorization of the same space $\CC^N$ by various subtori $T^K\subset T^N$ of the maximal torus, the $S^1$-equivariant Floer theories are related to each other in such a way that all the series $I^H_X$ can be expressed by complex oscillating integrals with {\em the same phase function} (``superpotetial''). This is how this was done in \cite{GiZ}.

Compute the function $I^H_X$ in the extreme case $X=\CC^N//T^N=pt$,
corresponding to the identity matrix $\m=(m_{ij})$, and $K=N$. We use $x_1,\dots, x_N$ in lieu of $Q_1,\dots, Q_N$, and take in account that $u_j=p_j=0$ in $H^*(pt)$: 
\[ I^H_{pt} = \sum_{D\in \ZZ^N} x^D\prod_{j=1}^N\frac{\prod_{r=-\infty}^0(-rz)}{\prod_{r=-\infty}^{D_j}(-rz)} =
\sum_{d\in \ZZ_{+}^N}\prod_{j=1}^N\frac{x_j^{D_j}}{(-z)^{D_j}D_j!}=e^{\textstyle-\sum x_j/z}.\]

{\tt Theorem 1} (see \cite{GiZ}). {\em The following complex oscillating integral
  \[ \I_X^H = \int_{\Gamma \subset {\small
      \left\{ \begin{array}{ll}\prod_j x_j^{m_{ij}}=Q_i \\ i=1,\dots,K \end{array} \right\} }} e^{\textstyle-\sum x_j/z} \, \frac{d\ln x_1\wedge \cdots \wedge d\ln x_N}{d\ln Q_1\wedge \cdots \wedge d\ln Q_K} \]
represents all solutions to the PDE system, $i=1,\dots, K$:
\begin{align*} \prod_{j: m_{ij}>0} & \prod_{r=0}^{m_{ij}-1}(rz-\sum_{i'=1}^Km_{i'j}zQ_{i'}\p_{Q_{i'}}) \ \I_X^H = \\
  Q_i \prod_{j: m_{ij}<0} & \prod_{r=0}^{-m_{ij}-1}(rz-\sum_{i'=1}^Km_{i'j}zQ_{i'}\p_{Q_{i'}})\ \I^H_X.\end{align*} }

{\tt Remarks.} (1) The system for $\I^H_X$ is written in the usual {\em scalar} representation of differential operators, while the system for $I^H_X$ was written in a {\em vector} representation, in which $-zQ_i\p_{Q_i}$ acts on functions with values in $H^*(X)$ by $p_i-zQ_i\p_{Q_i}$, where $p_i$ is the operator of multiplication in $H^*(X)$. Otherwise the differential systems are the same.

(2) The cycles $\Gamma$ here are non-compact {\em Lefschetz thimbles}, to be constructed by ascending gradient flow for the real part (for $z>0$) of the Morse function $\sum_j x_j$ under the constraints $\prod_j x_j^{m_{ij}}=Q_i$, $i=1,\dots,K$. The number of independent cycles is equal to the number of critical points and agrees with the rank of the differential system, but
exceeds, generally speaking, the dimension of $H^*(X)$. This is because the mirror integral and the differential system depends only on the subtorus $T^K\subset T^N$ but (unlike the series $I^H_X$), not on the choice of a chamber of regular values $\omega$ of the moment map. To select critical points and thimbles relevant for a given $X$, one needs to let ``quantum parameters'' $Q$ degenerate in such a way that $Q^d\to 0$ when $d$ lies in the Mori cone of $X$, and ignore those critical points which in this limit escape to infinity. This will become more clear when we compute the critical points.

(3) Similar mirror formulas can be written for toric bundle spaces $E\to X$ \cite{GiE}, super-bundles $\Pi E$ \cite{GiT, GiE}, and a certain class of toric complete intersections \cite{GiH, GiT}. We will touch upon these in the context of K-theory.

\section*{Toric mirrors in K-theory}

We retrace our steps in the K-theoretic context.

Take the toric $q$-hypergeometric function
\[  I^K_X=\sum_{d\in \ZZ^K} Q^d \prod_{j=1}^N \frac{\prod_{r=-\infty}^0(1-U_j(P)\gL_j^{-1}q^r)} {\prod_{r=\infty}^{D_j(d)}(1-U_j(P)\gL_j^{-1}q^r)}.\]
We recall from Part V that it takes values in the ring $K^0(X)$, which is multiplicatively generated by line bundles $P_1,\dots P_K$ with $\ch (P_i)=e^{-p_i}$, that $U_j$  are line bundles corresponding to the toric divisors $-u_j$, and that they satisfy multiplicative relations and Kirwan's relations
expressed respectively as
\[ U_j(P)=\prod_{i=1}^K P_i^{m_{ij}}, \ j=1,\dots,N, \ \ \text{and}\ \prod_{j\in J}(1-U_j)=0\ \text{if}\ \prod_{j\in J}u_j=0.\]
Now compute $I^K_X$ in the extreme case $X=\CC^N//T^N$:
\[ I^K_{pt}=\sum_{D\in \ZZ_{+}^N}\frac{X_1^{D_1}\cdots X_N^{D_N}}{\prod_{j=1}^N\prod_{r=1}^{D_j}(1-q^r)}=e^{\textstyle \sum_{k>0}\sum_{j=1}^N X_j^k/k(1-q^k)}.\]

{\tt Theorem 2.} {\em The following complex oscillating integral
\[ \I_X^K = \int_{\Gamma \subset {\small
      \left\{ \begin{array}{ll}\prod_j X_j^{m_{ij}}=Q_i \\ i=1,\dots,K \end{array} \right\} }}\  e^{\textstyle \sum_{k>0}\sum_j X_j^k/k(1-q^k)} \ \frac{\bigwedge_{j=1}^N d\ln X_j}{\bigwedge_{i=1}^K d\ln Q_i} \]
satisfies the system of $q$-difference equations, $i=1,\dots, K$:
\begin{align*} \prod_{j: m_{ij}>0} & \prod_{r=0}^{m_{ij}-1}(1-q^{-r}q^{\sum_{i'}m_{i'j}Q_{i'}\p_{Q_{i'}}}) \ \I_X^K = \\ Q_i \prod_{j: m_{ij}<0} & \prod_{r=0}^{-m_{ij}-1}(1-q^{-r}q^{\sum_{i'}m_{i'j}Q_{i'}\p_{Q_{i'}}})\ \I^K_X.\end{align*} }

{\tt Remark.} It is straightforward to check that the $K^0(X)$-valued $q$-hypergeometric series $I^K_X$ satisfies the same equations, except that the translation operators $q^{Q_i\p_{Q_i}}$ act in the vector representation by the composition $P_i q^{Q_i\p_{Q_i}}$ of translation and multiplication by $P_i$.

\medskip  

{\tt Proof.} We have
\[  q^{X\p_X} e^{\sum_{k>0} X^K/k(1-q^k)} =(1-X)e^{\sum_{k>0}X^k/k(1-q^k)},\]
and consequently for $r\geq 0$
\[ (1-q^{-r}q^{X\p_X})\, X^r e^{\sum_{k>0}X^k/k(1-q^k)} = X^{r+1} e^{\sum_{k>0}X^k/k(1-q^k)}.\]
Therefore for each $i=1,\dots, K$
\begin{align*} \prod_{m_{ij}>0}\prod_{r=0}^{m_{ij}-1}(1-q^{-r}q^{X_j\p_{X_j}})\, I^K_{pt}& = \prod_{m_{ij}>0}X_j^{m_{ij}}\ I^K_{pt}, \\
\prod_{m_{ij}>0}\prod_{r=0}^{-m_{ij}-1}(1-q^{-r}q^{X_j\p_{X_j}})\, I^K_{pt}& =Q_i \prod_{m_{ij}<0} X_j^{-m_{ij}}\ I^K_{pt},\end{align*}
which are the same under Batyrev constraints $\ln Q_i =\sum_j m_{ij} \ln X_j$.
Note that the volume form is translation-invariant, and that $X_j\p_{X_j}$ projects to $\sum_i m_{ij} Q_i \p_{Q_i}$. The result follows. \qed

\medskip

For the sake of completeness, let us include formulas for toric bundles an super-bundles. Let $V_a$, $a=1,\dots, L$ be {\em negative} line bundles, $V_a=\prod_{i=1}^K P_i^{l_{ia}}$, and $E=\oplus_{i=1}V_a$. Under these assumptions, the $q$-hypergeometric functions of Part V, namely $I_{E}$ and $I_{\Pi E^*}$ (note the dual) have well-defined non-equivariant limits:
\begin{align*} I_E&=\sum_{d\in \ZZ^K} Q^d \, \prod_{j=1}^N\frac{\prod_{r=-\infty}^0(1-q^rU_j(P))} {\prod_{r=-\infty}^{D_j(d)}(1-q^rU_j(P))}\ \prod_{a=1}^L\prod_{r=0}^{\Delta_a(d)-1}(1-q^{-r}V_a(P)) \\ 
  I_{\Pi E^*}&=\sum_{d\in \ZZ^K} Q^d \, \prod_{j=1}^N\frac{\prod_{r=-\infty}^0(1-q^rU_j(P))} {\prod_{r=-\infty}^{D_j(d)}(1-q^rU_j(P))}\  \prod_{a=1}^L \prod_{r=1}^{\Delta_a(d)}(1-q^r V_a(P)),\end{align*}
where $\Delta_a(d)=\sum_{i=1}^K d_il_{ia}$, and the sums effectively range over
$d$ with $D_j(d)\geq 0$ for each $j$. To save space, we write out the systems of $K$ $q$-difference equations (i.e. $i=1,\dots, K$), for $I_E$:
\[ \prod_j \prod_{r=0}^{m_{ij}-1}(1-q^{-r}q^{\sum_{i'}m_{i'j}Q_{i'}\p_{Q_{i'}}}) \I = Q_i\prod_a \prod_{r=0}^{l_{ia}-1}(1-q^{r}q^{-\sum_{i'}l_{i'a}Q_{i'}\p_{Q_{i'}}}) \I,\]
and for $I_{\Pi E^*}$:
\[ \prod_j\prod_{r=0}^{m_{ij}-1}(1-q^{-r}q^{\sum_{i'}m_{i'j}Q_{i'}\p_{Q_{i'}}}) \I  
= Q_i \prod_a \prod_{r=1}^{l_{ia}}(1-q^{-r}q^{\sum_{i'}l_{i'a}Q_{i'}\p_{Q_{i'}}}) \I,\]
following two conventions. Firstly, they are written in the scalar representation, i.e. $q^{Q_i\p_{Q_i}}$ acts on $I_E$ and $I_{\Pi E^*}$ by $P_iq^{Q_i\p_{Q_i}}$. Secondly, the products of the form $\prod_{r=0}^{m-1}$ should be interpreted as
$\prod_{r=-\infty}^{m-1}/\prod_{r=-\infty}^{-1}$, and when $m<0$, turn therefore into
$\prod_{r=m}^{-1}$ of {\em inverse} operators. To avoid inverses, one should move such factors to the other side of the equations, commuting them correctly across $Q_i$ if needed (as we did above in the formulations of the theorems).  

The K-theoretic mirrors for $E$ and $\Pi E^*$ are solutions to the above $q$-difference systems in the form of complex oscillating integrals
\begin{align*} \I_{E}&=\int_{\Gamma\subset \hat{X}_Q} 
e^{
  \sum_{k>0}(\sum_j X_j^k-q^k\sum_a Y_a^k)/k(1-q^k)} \ \text{\small $ \frac{\bigwedge_{j=1}^N d\ln X_j\, \bigwedge_{a=1}^Ld\ln Y_a}{\bigwedge_{i=1}^K d\ln Q_i}$}, \\
 \I_{\Pi E^*}&=\int_{\Gamma\subset \hat{X}_Q} 
e^{
  \sum_{k>0}(\sum_j X_j^k-q^k\sum_a Y_a^k)/k(1-q^k)} \ \text{\small $ \frac{\bigwedge_{j=1}^N d\ln X_j\, \bigwedge_{a=1}^LdY_a}{\bigwedge_{i=1}^K d\ln Q_i}$} \end{align*}
over Lefschetz' thimbles in the tori
\[ \hat{X}_Q=\left\{(X,Y)\in (\CC^{\times})^{N+L}\, \middle| \,  \prod_{j=1}^N X_j^{m_{ij}}=Q_i \prod_{a=1}^L Y_a^{l_{ia}} , i=1,\dots, K\right\}.\]

\section*{$D_q$-modules and quasiclassics}

A complex oscillating integral $\I(Q)=\int_{\gamma} e^{f_Q(x)/z} dv_Q$ is set by a phase function $f_Q$ and a volume form $dv_Q$, both depending on parameters $Q$. In the quasiclassical limit, i.e. as the wave length parameter $z$ tends to $0$,  the asymptotical expansion of such integrals has the form
\[ z^{n/2}\frac{e^{f_Q(x_{cr})/z}}{\sqrt{\Delta (x_{cr})}} \times \text{(power series in $z$)},\]
where $x_{cr}$ is a critical point of the phase function, $\Delta (x_{cr})$ is the Hessian of $f_Q$ with respect to the volume form at the critical point, and $n$ is the dimension of the integral. These geometric data are parameterized by critical points of the functions in the family. They naturally live on the Lagrangian variety $L\subset T^*B$ {\em generated} (in Arnold's terminology) by the family $f_Q$ in the cotangent bundle of the parameter space:
\[  L=\left\{ (p,Q) \in T^*B\ \middle| \ \exists x_{cr}:  p=d_Qf_Q(x_{cr}) \right\} .\]
Let us compute such Lagrangian varieties for Hori--Vafa mirrors of toric manifolds, first cohomological, and then K-theoretic.

We have the function $ \sum_j x_j$ under the constraints $\sum_j m_{ij}\ln x_j=\ln Q_i$, $i=1,\dots, K$. Introducing Lagrange multipliers $p_i$ and the auxiliary function of $(x,p)$:
\[ \sum_{j=1}^N x_j - \sum_{i=1}^K p_i \left(\sum_{j=1}^N m_{ij} \ln x_j-\ln Q_i\right) ,\]
and find the critical points 
\[ x_j = u_j(p) = \sum_{i=1}^K p_i m_{ij},\ \ j=1,\dots, K.\]
This describes $L$ as a Lagrangian variety with respect to the action $1$-form 
\[ p_1 \,\frac{dQ_1}{Q_1} + \cdots + p_K\, \frac{dQ_K}{Q_K},\]
given by {\em Batyrev relations}\footnote{Well, they are found in the 1992 preprint \cite{GiA}, where the multiplicative structure on Floer homology of toric manifolds $\CC^N//T^K$ was studied by means of discretization of loops in $\CC^N$.}\cite{Bat}:
\[ L=\left\{ (p,Q)\ \middle|\ Q_i=\prod_{j=1}^N u_j(p){m_{ij}} \right\} .\]
More precisely, as when all $Q_i\neq 0$, these are equivalent to the deformations of Kirwan's relations in the quantum cohomology algebra of $X$. 

Let us now turn to K-theory. Note that, though $q=e^z$ tends to $1$ as $z$ tends to $0$, the integral is also ``oscillating'' whenever $q$ tends to any root of unity. So, we examine the stationary phase asymptotics near $q=\zeta^{-1}$, where $\zeta$ is a primitive $m$-th root of unity. We single out the fast-oscillating part of the phase function in the integral:
\[ \int_{\Gamma\subset \hat{X}_Q} e^{\textstyle\frac{1}{1-q}\sum_{l>0} \sum_j \frac{X_j^{lm}}{l^2m^2}} \left( \text{Amplitude regular at $q=\zeta^{-1}$} \right) \times \left(\text{Volume}\right) . \]
The auxiliary function with Lagrange multipliers $p_i$ is
\[ \sum_{l>0} \sum_j \frac{X_j^{lm}}{l^2m^2} - \sum_i p_i\left( \sum_j m_{ij}\ln X_j - \ln Q_i\right) .\]
We compute the critical points: for $j=1,\dots, N$, 
\[ \sum_{l>0} \frac{X_j^{lm}}{lm}=\sum_i p_i m_{ij}=u_j(p), \ \ \text{or}\  \ln (1-X_j^m)^{1/m} = u_j(p).\]
Recalling the K-theoretic notation $P_i=e^{-p_i}$, $U_j(P)=\prod_i P_i^{m_{ij}}$,
we find that the critical points parameterize Lagrangian variety
\[ L_m = \left\{ (P, Q)\ \middle| \ Q_i^m = \prod_j (1-U_j^m(P))^{m_{ij}}, \ i=1,\dots, K \right\} \]
  in the complex torus with the symplectic structure
  \[ \frac{dP_1}{P_1}\wedge \frac{dQ_1}{Q_1}+\cdots+\frac{dP_K}{P_K}\wedge \frac{dQ_K}{Q_K}.\]
  There are several conclusions to draw.

(i) For $m=1$, the Lagrangian variety is given by ``quantum'' deformations of Kirwan's relation in $K^0(X)$.

(ii) Stationary phase asymptotics at different roots of unity exhibit self-similar behavior: $L_m$ is the inverse image of $L_1$ under the Adams
map $\Psi^m: (P,Q)\mapsto (P^m, Q^m)$.

(iii) Lagrangian varieties $L_m$ live in the symplectic torus (rather then the cotangent bundle space) which is a natural home for symbols of $q$-difference operators.

In view of these toric examples an their properties, three general problems arise: To develop a general theory of characteristic varieties of $D_q$-modules,
to describe the class of complex oscillating integrals which exhibit self-similar behavior when the spectral parameter approaches the roots of unity, and to include K-theory into the philosophy of mirror symmetry.

\enddocument